\documentclass[prl, reprint]{revtex4-1}

\usepackage[utf8]{inputenc}
\usepackage{amsmath,bm, romanbar, amsfonts, graphicx}
\usepackage{subfig}
\usepackage{booktabs}
\newcommand{\Qthree}{\mathcal{Q}_{\Romanbar{III}}}
\newcommand{\qthree}{q_{\Romanbar{III}}}
\newcommand{\hatQthree}{\hat{\mathcal{Q}}_{\Romanbar{III}}}
\newcommand{\Pdensity}{\Xi}

\begin{document}

\title{Family of closed-form solutions for two-dimensional correlated diffusion processes}

\author{Haozhe Shan}
\affiliation{Program in Neuroscience, Harvard University, Boston, Massachusetts, 02115, USA}
\affiliation{Center for Brain Science, Harvard University,  Cambridge, Massachusetts, 02138, USA}
\author{Rub{\'e}n Moreno-Bote}%
\affiliation{Department of Information and Communications Technologies, Pompeu Fabra University, Barcelona, Spain}
\affiliation{Center for Brain and Cognition, Pompeu Fabra University, Barcelona, Spain}
\affiliation{Serra H{\'u}nter Fellow Programme, Pompeu Fabra University, Barcelona, Spain}
\author{Jan Drugowitsch}
\affiliation{Department of Neurobiology, Harvard Medical School, Boston, Massachusetts, 02115, USA}
\affiliation{Center for Brain Science, Harvard University, Cambridge, Massachusetts, 02138, USA}
\email{jan\_drugowitsch@hms.harvard.edu}

\date{\today}

\begin{abstract}
Diffusion processes with boundaries are models of transport phenomena with wide applicability across many fields. These processes are described by their probability density functions (PDFs), which often obey Fokker-Planck equations (FPEs). While obtaining analytical solutions is often possible in the absence of boundaries, obtaining closed-form solutions to the FPE is more challenging once absorbing boundaries are present. As a result, analyses of these processes have largely relied on approximations or direct simulations. In this paper, we studied two-dimensional, time-homogeneous, spatially-correlated diffusion with linear, axis-aligned, absorbing boundaries. Our main result is the explicit construction of a full family of closed-form solutions for their PDFs using the method of images (MoI). We found that such solutions can be built if and only if the correlation coefficient $\rho$ between the two diffusing processes takes one of a numerable set of values. Using a geometric argument, we derived the complete set of $\rho$'s where such solutions can be found. Solvable $\rho$'s are given by $\rho =  - \cos \left( \frac{\pi}{k} \right)$, where $k \in \mathbb{Z}^+ \cup \{ +\infty\}$. Solutions were validated in simulations. Qualitative behaviors of the process appear to vary smoothly over $\rho$, allowing extrapolation from our solutions to cases with unsolvable $\rho$'s.
\end{abstract}

\maketitle

\section{Introduction}

Diffusion processes with absorbing boundaries are essential tools to model a multitude of real-world processes.
In neuroscience, for example, they act as models of decision-making \cite{bogacz2006physics} and neuronal action potential generation \cite{moreno2002response};
in finance, they are used for stock pricing \cite{dshalalow2005exit} and risk modeling \cite{yi2010first}; and in physics, they have, for example, been used to model movement of charges through conductors \cite{richert1989diffusion, hirao1995diffusion}.
Some applications involve higher-dimensional diffusions in which the process can become spatially correlated.
For example, decision-making models can assume multiple, correlated sources of decision-related evidence \cite{moreno2010decision}.
In models of neural action potentials, correlated diffusions might occur if neurons receive shared inputs \cite{moreno2006auto}.
In these cases, it is essential to understand how these correlations impact the process' behaviors. 

Due to their large number of applications, diffusion processes with boundaries have been widely studied.
A fundamental quantity describing these processes is the probability density function (PDF) of $\bm{x}(t)$, here denoted as $\Pdensity(\bm{x}, t)$.
Its time-evolution is described by the Fokker-Planck equation (FPE) \cite{coxmiller1965}.
Once the PDF is known by solving the FPE, multiple other process properties, such as the survival probability and boundary first-passage times, can be derived.
Therefore, finding the PDF is the first step towards a better understanding of the process.
Much work has been performed on finding solutions to the FPE for one-dimensional processes.
In this case, PDFs for problems with one or two time-invariant absorbing boundaries can be found with the method of images, which constructs PDFs as linear superpositions of free-space solutions \cite{coxmiller1965}.
In higher dimensions, however, the geometry of the process becomes significantly more complex, especially in the presence of a non-zero drift.
As a result, analytic results are scarce despite decades of effort.
Notable work includes \cite{iyengar1985hitting}, which provided analytic PDF expressions for two-dimensional processes with orthogonal boundaries on one side but without drift.
Later work provided some corrections, and proposed numerical methods to approximate the PDFs of processes with drift \cite{metzler2010first}.
Most relevantly for our work is \cite[Sec. 6]{sacerdote2016first}, which provided PDFs for drifting two-dimensional, correlated diffusion processes as infinite sums of Bessel functions with space- and time-dependent arguments. 

In particular for two-dimensional, correlated diffusion processes, the majority of previous work only provided analytical PDF expressions that gives limited insight, and whose numerical evaluation might be cumbersome.
Closed-form expressions are known only in rare cases, such as for uncorrelated processes or processes with a correlation coefficient of $\rho = -0.5$ \cite{moreno2010decision}. \cite{van2016common} derived the solution for $\rho=-\frac{\sqrt{2}}{2}$, but did not provide the resulting expression.
Such closed-form expressions have multiple benefits. They are usually easier to interpret than direct simulations, and they can provide a core around which analytical approximation can be found by perturbative expansions (e.g., colored noise diffusion from white noise \cite{moreno2002response, moreno2004role, moreno2006auto}).
For numerical analyses, closed-form expressions are significantly cheaper to compute than simulations and they are easy to evaluate to machine precision with finite operations.

Our aim was to find closed-form expressions for PDFs that describe two-dimensional, correlated diffusions with drift in the presence of  two time-invariant, orthogonal, and absorbing boundaries.
We approached this problem by focusing on solutions that can be constructed with the method of images (MoI). In particular, we aimed to determine under which circumstances we can find solutions that are expressible with such a finite number of images.
In what follows, we show that the only property that determines if such a solution exists is the diffusion process' correlation coefficient $\rho$.
Specifically, with the exception of $\rho = 1$, no close-form solution exists for positive $\rho$.
Furthermore, for negative $\rho$'s, we can only find solutions for a countable but infinite number of $\rho$.
For those, we provide the closed-form solutions, and demonstrate their validity in numerical simulations.
Thus, our work provides the complete set of all FPE solutions that can be found by the MoI and contain a countable number of images for this problem.

\section{Results}
We consider a 2D diffusion process with drift, denoted as  $\bm{x}(t)$. Its dynamics are given by
\begin{equation}
    d\bm{x}(t) =\bm{\mu}dt+\bm{\xi}(t),
    \label{eq:langevin}
\end{equation}
where $\bm{\mu}$ is the drift rate and $\bm{\xi}(t)$ is a Gaussian process.
$\bm{\xi}(t)$ has zero mean and covariance
\begin{equation}
    \left< \bm{\xi}(t) \bm{\xi}(t')^T\right>=\delta(t-t')\bm{\Sigma},
\end{equation}
where $\delta (\cdot)$ is the Dirac delta function and
\begin{equation}
    \bm{\Sigma} =\begin{pmatrix} 1 & \rho \\ \rho & 1 \end{pmatrix}  \quad \rho \in [-1,1].
\end{equation}
While we develop our solutions for covariance matrices of this specific form, our results also capture processes $\bm{\tilde{x}}(t)$ with arbitrary positive definite covariance matrices $\bm{\tilde{\Sigma}}$ by letting $x_1(t)\equiv \tilde{x_1}(t) / \sqrt{\tilde{\Sigma}_{11}}$ and $x_2(t)\equiv \tilde{x_2}(t) / \sqrt{\tilde{\Sigma}_{22}}$. The probability density function $\Pdensity(\bm{x},t)$ obeys the FPE
\begin{equation}
\frac{\partial \Pdensity}{\partial t}= - \bm{\mu} \cdot \nabla_{\bm{x}} \Pdensity + \frac{1}{2} \sum_{i,j=1}^2 \Sigma_{ij} \frac{\partial \Pdensity}{\partial x_i\partial x_j}.
\label{eq:fpe}
\end{equation}

We assume the initial condition $\bm{x}(0)=\bm{s}^{(0)}$. This is equivalent to
\begin{equation}
    \Pdensity(\bm{x},t=0)=\delta(\bm{x}-\bm{s}^{(0)}).
    \label{eq:initial}
\end{equation}
Without loss of generality, we assume the process to be in the third quadrant $\Qthree$ in Cartesian coordinates. Thus,
\begin{equation}
    \bm{s}^{(0)} \in \Qthree \equiv \left\{ \bm{x} | x_1 < 0, x_2 < 0 \right\} .
\end{equation}
The process is bounded from above by two linear, axis-aligned, absorbing boundaries at
\begin{align}
    B_1 &\equiv \left\{ \bm{x} | x_1 \leq 0, x_2 = 0 \right\} , \\
    B_2 &\equiv \left\{ \bm{x} | x_1 = 0, x_2 \leq 0 \right\},
\end{align}
such that it additionally needs to obey Dirichlet (a.k.a.~absorbing) boundary conditions
\begin{equation}    
    \forall t\geq0, \forall \bm{x} \in B_1 \cup B_2: \Pdensity(\bm{x},t)=0.
    \label{eq:boundary condition}
\end{equation}
These boundary conditions ensures that no probability mass enters the space outside the third quadrant (see \cite[Sec.~5.7, Eq.~(63)]{coxmiller1965}). Finally, the solution for $\Pdensity(\bm{x},t)$ must be non-negative everywhere. This is guaranteed by the boundary conditions and the maximum principle of elliptic PDEs, to which the considered FPE belongs \cite[Chap. 2]{han2011elliptic}. 

In what follows, we study the \textit{Dirichlet problem} of obtaining solutions ($\Pdensity(\bm{x},t)$) to Eq.~(\ref{eq:fpe}) under constraints Eqs.~(\ref{eq:initial}) and (\ref{eq:boundary condition}).
We will derive necessary and sufficient conditions for the existence of MoI solutions, and will determine these solutions in cases where they exist.
As we will show, such solutions only exist for a discrete set of correlation coefficients $\rho$.
To show this, we will first discuss a general expression for MoI constructions.
Second, we will identify the conditions under which MoI constructions satisfy the boundary condition, and, third, the initial condition.
Our approach reveals a restricted, discrete set of $\rho$'s for which exact solutions can be found, while at the same time providing closed-form expressions for these solutions.
Lastly, we will validate exact solutions with numerical simulations of the process.

\subsection{MoI construction of potential solutions}

Let us for now ignore the boundary condition, Eq.~(\ref{eq:boundary condition}).
In this case, free-space solutions to the FP equation are known to be the PDF of a bivariate Gaussian distribution with mean $\bm{s}^{(0)}+\bm{\mu}t$ and covariance $\bm{\Sigma}t$,
\begin{multline}
   \mathcal{N}\left(\bm{s}^{(0)}+\bm{\mu}t, \bm{\Sigma}t\right)  \\
   \equiv \frac{1} {2 \pi t \sqrt{| \bm{\Sigma} |}} e^{ - \frac{1}{2t} \left(\bm{x} - \bm{s}^{(0)} - \bm{\mu} t \right)^{T} \bm{\Lambda} \left( \bm{x} - \bm{s}^{(0)} - \bm{\mu} t \right)},
    \label{eq:image}
\end{multline}
where $\bm{\Lambda} = \bm{\Sigma}^{-1}$.

The MoI constructs solutions to a PDE with boundary conditions by adding scaled \textit{image functions}, $\mathcal{N}\left(\bm{s}^{(i)}+\bm{\mu}t, \bm{\Sigma}t\right)$, to Eq.~(\ref{eq:image}).  Hereafter, we will simply refer to the image function as an \textit{image} $\bm{s}^{(i)}$, and the point in space specified by $\bm{s}^{(i)}$ as a \textit{source} $\bm{s}^{(i)}$.
The resulting MoI construction has the form
\begin{multline}
    \Pdensity \left( \bm{x}, t; \bm{s}^{(0)} \right)
    = \mathcal{N}\left(\bm{s}^{(0)}+\bm{\mu}t, \bm{\Sigma}t\right) \\ + \sum_{i=1}^{N-1} a_i \mathcal{N}\left(\bm{s}^{(i)}+\bm{\mu}t, \bm{\Sigma}t\right),
    \label{eq:MoI construction}
\end{multline}
where $a_i$ is the image weight associated with image $\bm{s}^{(i)}$. 
MoI constructions of the form of Eq.~(\ref{eq:MoI construction}) satisfy the FP equation due to linearity of PDEs, and will satisfy Eq.~(\ref{eq:initial}) for any value of the weights if none of the added images are in the third quadrant. 
Furthermore, they will satisfy Eq.(\ref{eq:boundary condition}) if images cancel each other on the boundaries.
We denote  the set of sources by $\Omega=\{ \bm{s}^{(i)}\} _{i=0}^{N-1}$ and the set of image weights by $\kappa  = \left\{ a_{i} \right\}_{i=0}^{N-1}$ (with constant $a_0=1$). Overall, finding an expression of the form of Eq.~(\ref{eq:MoI construction}) that meets all criteria implies that we have identified a closed-form solution.

\subsection{Satisfying the boundary conditions}
\subsubsection{Placement of canceling images}

As discussed above, an MoI construction satisfies the boundary condition if and only if images cancel each other at the boundaries. More specifically, for any image with source $\bm{s}^{(i)}\in\Omega$ and at any time $t$, there should be a set of images that cancel it on both boundaries ($\{ B_{1,2}\}$).

Importantly, linearity of exponential functions in Eq.~(\ref{eq:MoI construction}) requires all images cancelling each other on a particular boundary to have the same exponent on that boundary for all times $t \ge 0$.
For instance, to cancel image $\bm{s}^{(i)}$ on $B_1$, we must add at least one other image $\bm{s}^{(j)}$ to $\Omega$ that satisfies
\begin{multline}
    \left(\bm{x} - \bm{s}^{(i)} - \bm{\mu} t \right)^{T} \bm{\Lambda} \left( \bm{x} - \bm{s}^{(i)} - \bm{\mu} t \right)\\
    = \left(\bm{x} - \bm{s}^{(j)} - \bm{\mu} t \right)^{T} \bm{\Lambda} \left( \bm{x} - \bm{s}^{(j)} - \bm{\mu} t \right) ,
\end{multline}
for all $t \ge 0$ and $\bm{x} \in B_1$.

\begin{figure}[ht]
    \centering
    \includegraphics[width=0.24\textwidth]{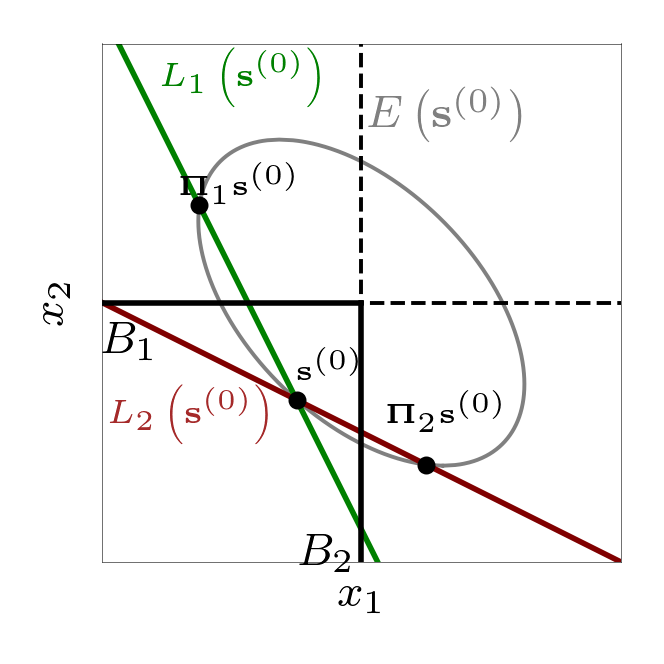}
    \caption{Example construction to cancel $\bm{s}^{(0)}$ at the two boundaries, here for $\rho = - \frac{1}{2}$ and some arbitrary $\bm{s}^{(0)}$.
        The intersections of $L_{1,2}$ with $E$ uniquely determine the placement of image sources to cancel $\bm{s}^{(0)}$ at $B_{1,2}$.}
    \label{fig:line_ellipse}
\end{figure}

To determine canceling images that satisfy this condition, we rewrite Eq.~(\ref{eq:image}) by separating terms linear and quadratic in $\bm{s}$ in the exponent,
\begin{multline}
    \mathcal{N}\left(\bm{s}+\bm{\mu}t, \bm{\Sigma}t\right) = \frac{1}{2 \pi t \sqrt{| \bm{\Sigma} |}} e^{ -\frac{1}{2t} \left( \bm{x} - \bm{\mu} t \right)^{T} \bm{\Lambda} \left( \bm{x} - \bm{\mu} t \right)} \\
    \times e^{- \frac{1}{2t} \bm{s}^T \bm{\Lambda} \bm{s}} \times e^{\frac{1}{t} \left(\bm{x} - \bm{\mu} t \right)^{T} \bm{\Lambda} \bm{s}}.
\end{multline}
The first exponent is independent of $\bm{s}$, and therefore shared by all images; matching the second exponent requires 
\begin{equation}
    \bm{s}^{(j)} \in E \left( \bm{s}^{(i)} \right) = \left\{ \bm{x} | \bm{x}^{T} \bm{\Lambda} \bm{x} = {\bm{s}^{(i)}}^{T} \bm{\Lambda} \bm{s}^{(i)} \right\};
\end{equation}
 matching the last exponent for $\bm{x}\in B_{1,2}$ respectively requires
 \begin{align}
    \bm{s}^{(j)} \in L_{1,2} \left( \bm{s}^{(i)} \right) = \left\{ \bm{x} \Big| \bm{e}_{1,2}^T \bm{\Lambda} \bm{x} = \bm{e}_{1,2}^T \bm{\Lambda} \bm{s}^{(i)} \right\},
\end{align}
where $\bm{e}_1 = (1, 0)^T$ and $\bm{e}_2 = (0, 1)^T$ are the Cartesian basis vectors.
Geometrically, $E \left( \bm{s}^{(i)} \right)$ is an ellipse and $L_{1,2} \left( \bm{s}^{(i)} \right)$ are lines, all of which pass through $\bm{s}^{(i)}$ (Fig.\ref{fig:line_ellipse}). The ellipse and each of the lines intersect at $\bm{s}^{(i)}$ and, in general, another point $\bm{s}^{(j)}$, leading to a unique canceling image.

Algebraically, it is easy to show that the mapping from to-be-cancelled image source $\bm{s}^{(i)}$ to canceling image source $\bm{s}^{(j)}$ for boundaries $B_{1,2}$, respectively, is given by
\begin{align}
    B_1: \quad \bm{s}^{(j)} &= \begin{pmatrix}
        1 & -2\rho \\
        0 & -1
    \end{pmatrix} \bm{s}^{(i)} \equiv \bm{\Pi}_1 \left(\rho\right)  \bm{s}^{(i)}, 
    \label{eq:pi_plus}\\
    B_2: \quad \bm{s}^{(j)} &= \begin{pmatrix}
        -1 & 0 \\
        -2 \rho & 1
    \end{pmatrix} \bm{s}^{(i)} \equiv \bm{\Pi}_2 \left(\rho\right)  \bm{s}^{(i)}.
    \label{eq:pi_minus}
\end{align}
Hereafter, we drop the dependency of the mapping $\bm{\Pi}_{1,2}$ on $\rho$ for notational convenience.
Both $\bm{\Pi}_{1,2}$ are involutory, that is $\bm{\Pi}_k^{-1} = \bm{\Pi}_k$ for both $k \in \{1, 2\}$.

To find the scaling coefficient $a_j$ for the canceling image, we solve
\begin{equation}
    a_i \mathcal{N}\left(\bm{s}^{(i)}+\bm{\mu}t, \bm{\Sigma}t\right) = -a_j \mathcal{N}\left(\bm{s}^{(j)}+\bm{\mu}t, \bm{\Sigma}t\right) ,
\end{equation}
for all $t \ge 0$ and for all $\bm{x} \in B_1$ or $\bm{x} \in B_2$, leading to
\begin{equation}
 a_j= - a_i e^{\bm{\mu}^T \bm{\Lambda} \left( \bm{s}^{(j)} - \bm{s}^{(i)}\right)} .
    \label{eq:scaling req}
\end{equation}

What would happen if $E \left( \bm{s}^{(i)} \right)$ and $L_k \left( \bm{s}^{(i)} \right)$ (for $k \in \{1, 2\}$) only intersect at a single point?
In this case, the line $L_k \left( \bm{s}^{(i)} \right)$ would be a tangent to the ellipse at $\bm{s}^{(i)}$.
It is easy to show that, for $L_1$ (or $L_2$), this only occurs if $s^{(i)}_2 = 0$ (or $s^{(i)}_1 = 0$), that is, if the image to be canceled happens to be located on one of the axes.
In those cases, $\bm{s}^{(j)} = \bm{\Pi}_k \bm{s}^{(i)} = \bm{s}^{(i)}$, such that the canceling image is mapped onto the image to be cancelled.
Furthermore, they receive opposite weights, that is $a_j = - a_i$, such that they cancel each other, and could be both removed.
While this is an intuitively odd scenario, it does not invalidate our approach, as the mappings $\bm{\Pi}_{1,2}$ remain valid.
Furthermore, as will become apparent later, no valid solution will have this property.
Therefore, it isn't a case that requires special attention.

\subsubsection{Finding a complete set of images}

\begin{figure}[ht]
    \centering
    \includegraphics[width=0.48\textwidth]{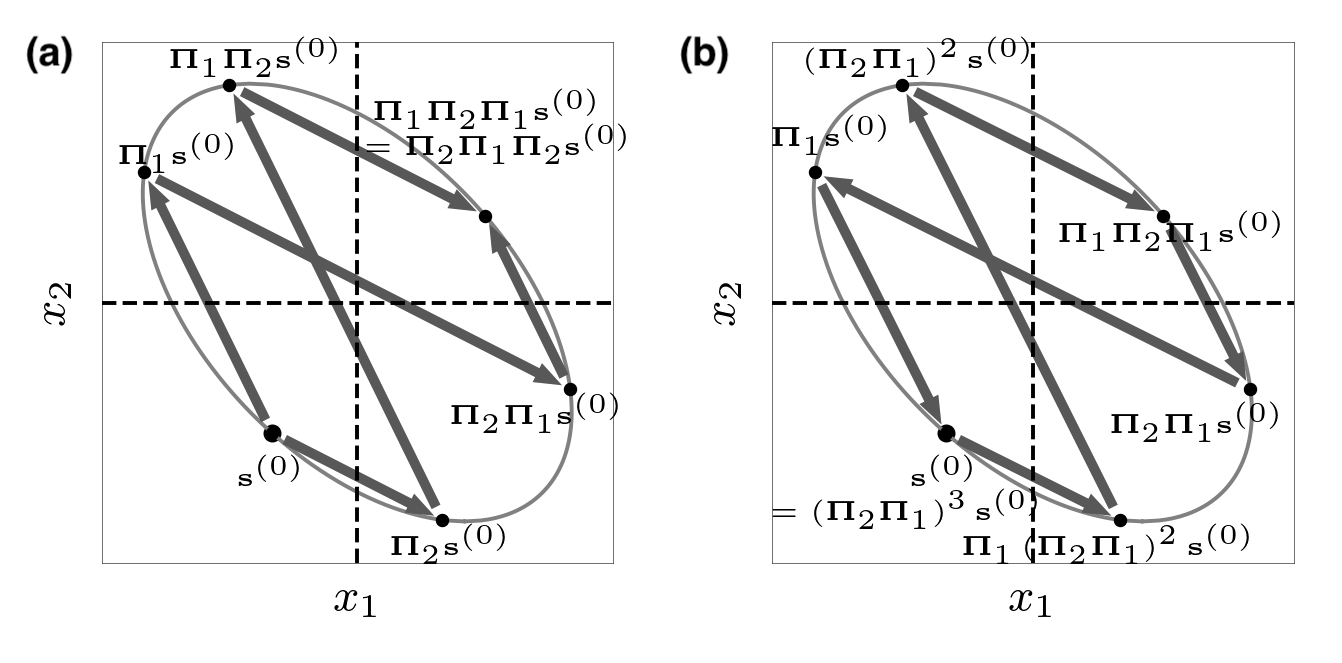}
    \caption{Two constructions to create a complete set of images, here illustrated for $\rho=-\frac{1}{2}$.
        Note that the same source can be expressed in multiple ways, allowing two formalisms to characterize the same set of images. (a) Pair-wise construction. (b) Sequential construction.}
    \label{fig2}
\end{figure}

We now consider how a set of images following the MoI construction can satisfy the boundary condition in Eq.~(\ref{eq:boundary condition}).
As we have two boundaries, two additional images, $\bm{\Pi}_1 \bm{s}^{(0)}$ and $\bm{\Pi}_2 \bm{s}^{(0)}$, are introduced to cancel the density from image $\bm{s}^{(0)}$.
The density of image $\bm{\Pi}_1 \bm{s}^{(0)}$ is canceled on boundary $B_1$ by image $\bm{s}^{(0)}$; however, it introduces additional density at $B_2$. Similarly, image $\bm{\Pi}_2 \bm{s}^{(0)}$ introduces some additional density at $B_1$. Therefore, yet another pair of images, $\bm{\Pi}_2 \bm{\Pi}_1 \bm{s}^{(0)}$ and $\bm{\Pi}_1 \bm{\Pi}_2 \bm{s}^{(0)}$, are needed, and a further pair of images to cancel their densities, and so on.
This leads to a pair-wise construction of images (Fig.~\ref{fig2}(a)). 

We proceed with the ansatz that having a finite set of images, that is $|\Omega|<\infty$, is a necessary condition for the existence of solutions (except when $\rho=-1$).
The motivating intuition is that one cannot introduce infinitely many sources without placing any in the third quadrant, thus violating the initial condition.
In a later section, this is shown to indeed be the case.

The number of images is finite if and only if, at some point, new images to be added are already in the set. Following the pair-wise construction (Fig.~\ref{fig2}(a)), this occurs only if there exist images that cancel two other images across different boundaries, that is
 \begin{equation}
     \exists \bm{s}^{(i,j,k)}\in\Omega:\bm{s}^{(k)} = \bm{\Pi}_1 \bm{s}^{(i)} = \bm{\Pi}_2 \bm{s}^{(j)} .
     \label{eq:meeting of pair}
 \end{equation}
The alternative of canceling two images across the same boundary, that is $\bm{\Pi}_1 \bm{s}^{(i)} = \bm{\Pi}_1 \bm{s}^{(j)}$, is invalid because this would imply $\bm{s}^{(i)}=\bm{s}^{(j)}$, in which case additional images are not necessary.

Due to the involutoriness of $\bm{\Pi}_{1,2}$, Eq.~(\ref{eq:meeting of pair}) leads to $\bm{s}^{(i)} = \bm{\Pi}_1 \bm{s}^{(k)}$ and $\bm{s}^{(j)} = \bm{\Pi}_2 \bm{s}^{(k)}$.
Here, $\bm{s}^{(i)}$ is the image that is cancelled by image  $\bm{s}^{(k)}$ across $B_1$.
Unless $\bm{s}^{(i)}$ is the original image $\bm{s}^{(0)}$, $\bm{s}^{(i)}$ was introduced to cancel another image $\bm{\Pi}_2 \bm{s}^{(i)} = \bm{\Pi}_2 \bm{\Pi}_1 \bm{s}^{(k)}$ across $B_2$.
This implies that, by further alternating application of $\bm{\Pi}_1$ and $\bm{\Pi}_2$, we can trace its origin back to $\bm{s}^{(0)}$.
It follows that if two images "meet" in the fashion of Eq.~(\ref{eq:meeting of pair}), then the complete set of images can be generated by following the sequence $\bm{\Pi}_1 \bm{s}^{(0)}, \bm{\Pi}_2 \bm{\Pi}_1 \bm{s}^{(0)}, \bm{\Pi}_1 \bm{\Pi}_2 \bm{\Pi}_1 \bm{s}^{(0)}, \dots$ until we find an image with source that coincides with the original location $\bm{s}^{(0)}$ (Fig.~\ref{fig2}(b)), resulting in $|\Omega|<\infty$.
Formally, this sequential image construction is given by the \textit{source generating function}
\begin{equation}
    \bm{\Phi} \left(\bm{s}^{(0)}, n \right) = 
    \begin{cases}
        \left(\bm{\Pi}_2 \bm{\Pi}_1\right)^{n/2} \bm{s}^{(0)} & \text{if } n \text{ is even}, \\
        \bm{\Pi}_1 \left(\bm{\Pi}_2 \bm{\Pi}_1\right)^{(n-1)/2} \bm{s}^{(0)} & \text{if } n \text{ is odd} .\\
    \end{cases}
    \label{eq:image gen}
\end{equation}
Following this formalism, we define all images with even $n$ "even-numbered images" (and similarly for "odd-numbered images").
For the number of images to be finite, we require
\begin{equation}
    \exists n \in \mathbb{Z}^+:\bm{\Phi} \left(\bm{s}^{(0)}, n \right) = \bm{s}^{(0)}.
    \label{eq:return to s0}
\end{equation}
This condition cannot hold for odd $n$, as $| \mathbf{\Pi}_j | =-1$ for both $j \in \{1, 2\}$, such that the product of an odd number of these mapping has determinant of -1, which cannot equate identity.
Therefore, we consider only even $n$.
Letting $n=2k$, Eq.~(\ref{eq:return to s0}) equals
\begin{equation}
    \exists k\in \mathbb{Z}^+:\left(\bm{\Pi}_2 \bm{\Pi}_1 \right)^k = \mathbb{I}.
    \label{eq:rotation}
\end{equation}
To show which values of $\rho$ satisfy this condition, let us first introduce a whitened process for mathematical convenience.

\subsubsection{Whitening the process}

We have shown that all images lie on the ellipse $E \left( \bm{s}^{(0)} \right)$ whose eccentricity increases with $|\rho|$.
We can simplify the analysis of image locations by whitening the process, after which all images come to lie on a circle, such that the location of each image is fully determined by its angle.
To perform this whitening, we desire to find a linear mapping $\bm{Q}$ that maps the original process $\bm{x}(t)$ into its whitened equivalent, $\hat{\bm{x}}(t) = \bm{Q} \bm{x}(t)$.
Under the required whitening constraint that $\bm{\Lambda} = \bm{Q}^T \bm{Q}$, one choice of $\bm{Q}$ that is symmetric is given by
\begin{equation}
 \bm{Q} = q \begin{pmatrix}
        \rho & \sqrt{1 - \rho^2} - 1 \\
        \sqrt{1 - \rho^2} - 1 & \rho 
    \end{pmatrix} , \label{eq:whitening matrix}\\
\end{equation}
with
\begin{equation}
 q = \frac{\text{sgn}(\rho)}{\sqrt{2 (1 - \rho^2) \left(1 - \sqrt{1 - \rho^2} \right)}} .
\end{equation}
This mapping, $\bm{Q}$, has singularities at $\rho=\pm 1$, where the process collapses into a 1D process with one or two absorbing boundaries for $\rho=1$ and $\rho=-1$, respectively.
Solutions in these special cases are known, and given in \cite[Sec.~5.7, Eqs.~(71) and (78)]{coxmiller1965}.
We thus restrict our discussions to $ 0<\lvert \rho \rvert<1$.

\begin{figure}[ht]
    \centering
    \includegraphics[width=0.48\textwidth]{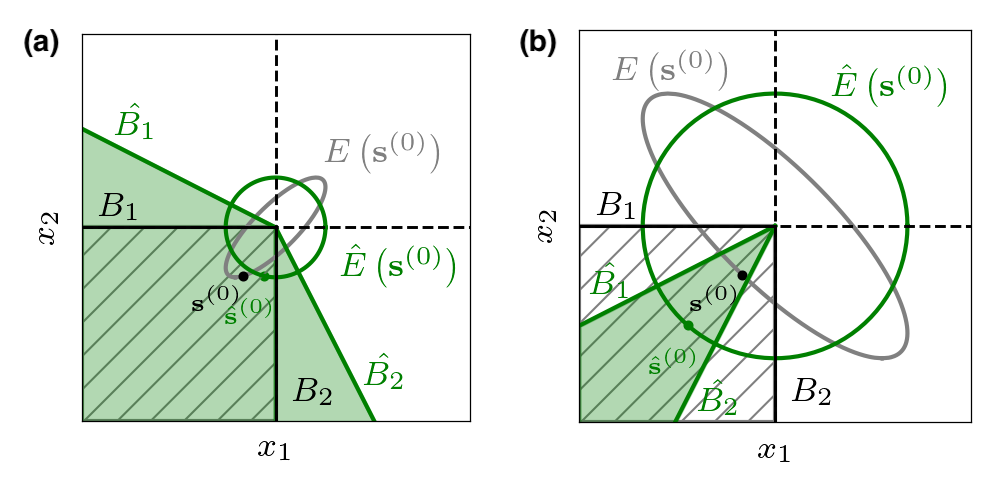}
    \caption{Whitening the diffusion (a) expands (for $0<\rho<1$)  or (b) shrinks (for $-1<\rho<0$) the original third quadrant (shaded) into a new region (green). The ellipse (gray line) on which images are found becomes a circle (green line). Boundaries ($B_1,B_2$, thick black line) are rotated into or away from the third quadrant ($\hat{B_1},\hat{B_2}$, thick green line). $\bm{s}^{(0)}$ here was chosen to not lie on the identity line, to show the resulting angular displacement of $\hat{\bm{s}}^{(0)}$ relative to $\bm{s}^{(0)}$ in such circumstances.}
    \label{fig:whitening}
\end{figure}

The consequences of this re-mapping are as follows (see Appendix for derivations).
First, the image cancellation maps for the whitened process become $\hat{\bm{\Pi}}_{1,2}=\bm{Q}\bm{\Pi}_{1,2}\bm{Q}^{-1}$, that, as before, obey $\hat{\bm{\Pi}}_1^2 = \hat{\bm{\Pi}}_2^2 = \mathbb{I}$.
Second, as desired, all images are now located on the circle
\begin{equation}
    \hat{\bm{x}}^T \hat{\bm{x}} = \hat{\bm{s}}^{(0) T} \hat{\bm{s}}^{(0)} ,
\end{equation}
where $\hat{\bm{s}}^{(0)} = \bm{Q} \bm{s}^{(0)}$ (Fig.~\ref{fig:whitening}).
Third, $\hat{\bm{\Pi}}_2 \hat{\bm{\Pi}}_1$ is a clockwise rotation matrix of angle $2\alpha$ (which we denote as  $\bm{R}(2\alpha)$), where $\alpha = \arccos(\rho) \in (0,\pi)$).
Thus, if we define the source generating function for the whitened process analogous to Eq.~(\ref{eq:image gen}) by
\begin{equation}
    \hat{\bm{\Phi}} \left(\hat{\bm{s}}^{(0)}, n \right) = 
    \begin{cases}
        \left(\hat{\bm{\Pi}}_2 \hat{\bm{\Pi}}_1\right)^{n/2} \hat{\bm{s}}^{(0)} & \text{if } n \text{ is even}, \\
        \hat{\bm{\Pi}}_1 \left(\hat{\bm{\Pi}}_2 \hat{\bm{\Pi}}_1\right)^{(n-1)/2} \hat{\bm{s}}^{(0)} & \text{if } n \text{ is odd},\\
    \end{cases}
    \label{eq:whitened image construction}
\end{equation}
then all even-numbered sources satisfy
\begin{equation}
    \hat{\bm{s}}^{(2m)} = \left(\hat{\bm{\Pi}}_2 \hat{\bm{\Pi}}_1\right)^m \hat{\bm{s}}^{(0)} = \bm{R}(2 m \alpha)   \hat{\bm{s}}^{(0)}
    \label{eq:whitened_even_images}
\end{equation}
for positive integers $m \in \mathbb{Z}^+$.
Fourth, the odd-numbered sources can be found similarly by
\begin{equation}
    \hat{\bm{s}}^{(2m + 1)} = \hat{\bm{\Pi}}_1 \hat{\bm{s}}^{(2m)} = \bm{F} \left[ \bm{R} \left( (2m + 1) \alpha \right) \hat{\bm{s}}^{(0)} \right] ,
    \label{eq:whitened_odd_images}
\end{equation}
which corresponds to a clockwise rotation of $\hat{\bm{s}}^{(0)}$ by $(2m + 1) \alpha$, followed by a flip across the anti-diagonal,
\begin{equation}
    \bm{F}= \begin{pmatrix}
        0 & - 1 \\
        - 1 & 0 
    \end{pmatrix}.
    \label{eq:antidiag flip}
\end{equation}

With these properties established, let us return to the question about which values of $\rho$ lead to a finite set of images, $| \Omega | < \infty$.
As for the non-whitened process, this number is finite if the source generating process returns to its origin after a finite number of steps.
We have already established that this only holds for an even number of steps.
Therefore, there needs to exists some integer $k$ such that $\hat{\bm{s}}^{(2k)} = \hat{\bm{s}}^{(0)}$.
By Eq.~(\ref{eq:whitened_even_images}), this implies $\bm{R}(2 k \alpha) = \mathbb{I}$, which holds as long as $2 k \alpha$ is some multiple of $2 \pi$ (i.e., one or several full rotations).
Overall, this means that the only values of $\rho$ that lead to a finite set of images, are
\begin{equation}
    \rho =  \cos \left( \alpha \right), \quad
    \textrm{with } \alpha =  \frac{l \pi}{k}  \quad k\in\mathbb{Z}^+,l = 1,\dots,k-1 ,
    \label{eq:candidate rho}
\end{equation}
where we have used $\alpha = \arccos(\rho)$, and have restricted $l$ to $1 \le l \le k-1$ to ensure $l / k < 1$.

\subsection{Satisfying the initial condition}

So far we have focused on satisfying the boundary condition, Eq.~(\ref{eq:boundary condition}), which has led to a restrictions on the values that $\rho$ can take.
Let us now consider which of those $\rho$'s additionally satisfy the initial condition, Eq.~(\ref{eq:initial}).
This condition implies that no other image than $\bm{s}^{(0)}$ can lie in the third quadrant.
After whitening, the third quadrant $\Qthree$ is mapped into a circular sector
\begin{multline}
    \hatQthree = \Bigg\{ \hat{\bm{x}} \Big| \hat{x}_2 \le \frac{\sqrt{1 - \rho^2} - 1}{\rho} \hat{x}_1 , \\
    \text{sgn}(\rho)\hat{x}_1 \le \text{sgn}(\rho)\frac{\sqrt{1 - \rho^2} - 1}{\rho} \hat{x}_2 \Bigg\} .
\end{multline}
For $\rho < 0$ or $\rho > 0$, the corresponding boundaries are rotated into or out of the original third quadrant, respectively (Fig.~\ref{fig:whitening}), by
\begin{equation}
    \label{eq:q3rot}
    \psi = \arctan \left( \frac{\sqrt{1 - \rho^2} - 1}{\rho} \right) = \frac{\alpha}{2} - \frac{\pi}{4} .
\end{equation}
Therefore, the angular width of $\hatQthree$ is $\pi/2 - 2 \psi = \pi - \alpha$.

In order to determine when it is possible to avoid placing sources (other than $\hat{\bm{s}}^{(0)}$) in the whitened third quadrant, $\hatQthree$, we will use the fact that all sources in the whitened space are located on a circle, such that it is sufficient to describe any source $\hat{\bm{s}}(i)$ by its polar angle $\theta_i$.
In addition, we will use the polar representation for the third quadrant,
\begin{equation}
    \hatQthree = \left\{ (r, \theta) \Big| \theta \in q_{\Romanbar{III}} \right\}  \textrm{ with }
    \qthree = \left[ \frac{3 \pi}{4} + \frac{\alpha}{2}, \frac{7 \pi}{4} - \frac{\alpha}{2} \right], 
    \label{eq:qthree_angle}
\end{equation}
where $\qthree$ is the range of polar angles within the third quadrant in whitened space.
A source $\hat{\bm{s}}(i)$ falls into the third quadrant if $\theta_i \in \qthree$.

\subsubsection{An infinite number of images}

We can now revisit the previous ansatz that, unless $\rho = -1$, sets with an infinite number of images will violate the initial condition.
To do so, note that, by Eqs.~(\ref{eq:whitened_even_images}) and (\ref{eq:candidate rho}), consecutive even-numbered images are placed at angular distance $2  \alpha = \frac{2 l \pi}{k}$ of each other.
However, as a complete construction might imply multiple full rotations (if $l > 1$), these consecutive even-numbered sources are not necessarily the even-numbered sources closest to each other.
Indeed, by periodicity, the angular spacing between all even-numbered sources is regular (Fig.~\ref{fig:imageangles}) and given by
\begin{equation}
    \beta = \frac{2 \pi}{k} .
\end{equation}
Therefore, once the number of images approaches infinity, $\beta$ approaches zero.
However, for any $| \rho | < 1$, $\alpha < \pi$ such that the angular width of $\hatQthree$ remains positive.
This implies that some images will fall into the third quadrant, such that such an infinite number of images will violate the initial condition.
This argument does not apply to $\rho = -1$, for which the whitening transformation is not well-defined.

\subsubsection{A finite number of images}

\begin{figure}[ht]
    \centering
    \includegraphics[width=0.48\textwidth]{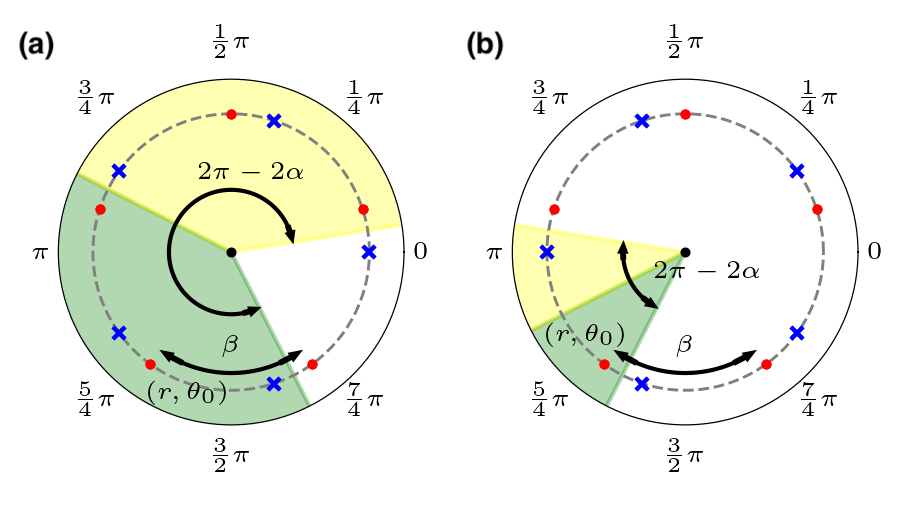}
    \caption{Using polar coordinates, the third quadrant after whitening, $q_{\Romanbar{III}}$ (green shade), and region D (yellow shade) together have width $2\pi-2\alpha$; geometrically adjacent even-numbered sources (red dots) have angular distance $\beta$. The odd-numbered sources are indicated by blue crosses. (a) $0<\rho<1$ ($\alpha=\frac{1}{5}\pi$). (b) $-1<\rho<0$ ($\alpha=\frac{4}{5}\pi$).}
    \label{fig:imageangles}
\end{figure}

Let us now focus on finite image sets.
Since, by Eq.~(\ref{eq:whitened_odd_images}), there is a one-to-one mapping from source $\hat{\bm{s}}^{(2m)}$ to source $\hat{\bm{s}}^{(2m+1)}$, we can reformulate the constraint that $\hat{\bm{s}}^{(2m+1)}\notin \hatQthree$ as a constraint on source $\hat{\bm{s}}^{(2m)}$.
That is, there exists some region $D$ such that $\theta_{2m} \notin D$ guarantees that $\theta_{2m + 1} \notin \qthree$. This region is given by (see Appendix)
\begin{equation}
    D = \left[ -\frac{ \pi}{4} + \frac{3 \alpha}{2},  \frac{\alpha}{2} + \frac{3\pi}{4} \right].
    \label{eq:avoidance region}
\end{equation}
Therefore, we only need to make sure that the even-numbered sources do not fall into $q_{\Romanbar{III}} \cup D$.
Incidentally, $q_{\Romanbar{III}}$ and $D$ are adjacent, such that they together from a single sector with combined width $2\pi- 2\alpha$ (Fig.~\ref{fig:imageangles}).
For what follows we will again exclude the special cases of $\rho\in\{\pm1, 0\}$ where solutions are known. This corresponds to taking $l\in{1,...,k-1}$ in Eq.~(\ref{eq:candidate rho}).

\begin{table*}
  \centering
  \begin{tabular}{@{}c@{\hspace{1cm}}c@{}} \toprule
  image mapping formalism & image rotation formalism \\ \midrule
    \parbox[][][t]{0.3\textwidth}{
        \begin{align*}
            \bm{s}^{(j)} &= \begin{cases}
                \left(\bm{\Pi}_2 \bm{\Pi}_1\right)^{j/2} \bm{s}^{(0)} & j \text{ is even},\\
                \bm{\Pi}_1 \left(\bm{\Pi}_2 \bm{\Pi}_1\right)^{(j-1)/2} \bm{s}^{(0)} & j \text{ is odd},\\
            \end{cases} \\
       \bm{\Pi}_1 &= \begin{pmatrix}
        1 & -2\rho \\
        0 & -1
    \end{pmatrix}, \quad \bm{\Pi_2} = \begin{pmatrix}
        -1 & 0 \\
        -2 \rho & 1
    \end{pmatrix} ,
        \end{align*}
    }
    &
    \parbox[][][t]{0.55\textwidth}{
        \begin{align*}
            \bm{s}^{(j)} &= \frac{1}{\sin \left( \frac{\pi}{k} \right)} \begin{cases}
             \begin{pmatrix}
        \sin \left( j \alpha + \frac{\pi}{k} \right) & \sin \left( j \alpha \right) \\
        -\sin \left( j \alpha \right) & - \sin \left( j \alpha - \frac{\pi}{k} \right)
    \end{pmatrix} \bm{s}^{(0)} & j \text{ is even},\\
    ~\\
          \begin{pmatrix}
        \sin \left( j \alpha \right) & \sin \left( j \alpha -  \frac{\pi}{k} \right) \\
        - \sin \left( j \alpha + \frac{\pi}{k} \right) & -\sin \left( j \alpha \right)
    \end{pmatrix} \bm{s}^{(0)} & j \text{ is odd},
            \end{cases} \\
            \alpha &= \frac{k-1}{k} \pi ,
        \end{align*}} \\
        \multicolumn{2}{c}{\parbox{0.85\textwidth}{
        \begin{align*}
                   \Pdensity(\bm{x},t)&=\mathcal{N}\left(\bm{s}^{(0)}+\bm{\mu}t, \bm{\Sigma}t\right)+\sum_{j=1}^{2k-1} a_j \mathcal{N}\left(\bm{s}^{(j)}+\bm{\mu}t, \bm{\Sigma}t\right),\\
            a^{(j)} &= \left(-1\right)^{j} \exp \left(\bm{\mu}^T \bm{\Lambda} \left( \bm{s}^{(j)} - \bm{s}^{(0)}\right) \right),\quad \textrm{where } \bm{\Lambda} = \bm{\Sigma}^{-1}.
        \end{align*}
    }}
  \end{tabular}
  \caption{Full expression of closed-form solutions for $2k$ images, corresponding to correlation coefficient $\rho = -\cos \left( \pi / k \right)$.  Image locations are given in two alternative formalisms. The image mapping formalism is based on Eq.~(\ref{eq:image gen}). The image rotation formalism is based on Eqs.~(\ref{eq:whitened_even_images}) and (\ref{eq:whitened_odd_images}) (see Appendix for derivation). Both solutions have a computational complexity that scales linearly with the number of images.}
  \label{tbl:solution_summary}
\end{table*}

As the angular width of both $q_{\Romanbar{III}}$ and $D$ individually is $\pi - \alpha$, a necessary condition for even-numbered sources to "avoid" them is for the angle between two geometrically adjacent sources $\beta$ to be larger than this region, that is
\begin{equation}
    \beta=\frac{2\pi}{k} > \pi - \alpha \qquad
    \Rightarrow \qquad k-l<2
\end{equation}
where the second inequality follows from $\alpha=l\pi/k$.
We will consider the two cases of positive and negative correlations $\rho$ in turn.

For all $0<\rho<1$, we have the additional constraint of $2l<k$ (since $\rho=\cos \left( \frac{l\pi}{k} \right)$).
Together with the previous inequality, this implies $l<2$ and thus $l=1$.
However, there is no value of $k\in\mathbb{Z}^+$ that satisfies both inequalities.
Thus, initial conditions cannot be satisfied for any  $0<\rho<1$, \textit{regardless of} $\hat{\bm{s}}^{(0)}$.

Since the original source is in the sector $q_{\Romanbar{III}}$, a sufficient condition for all other sources to skip the sector is
\begin{equation}
    \beta=\frac{2\pi}{k} \geq 2\pi-2\alpha \qquad
    \Rightarrow \qquad k-l\leq 1
\end{equation} 
which holds if we simply choose $l=k-1$. On the other hand, a necessary condition is that $\beta$ is larger than half the width of the combined $q_{\Romanbar{III}} \cup D$ region, such that even if $\hat{\bm{s}}^{(0)}$ approaches the center of that region, adjacent sources will be outside of it.
This leads to
\begin{equation}
    \beta=\frac{2\pi}{k} > \pi-\alpha \qquad \Rightarrow \qquad k-l< 2,
\end{equation} 
which only holds for $l=k-1$.
Thus, for $-1<\rho<0$, initial conditions are satisfied if and only if $l=k-1$.

\subsection{Exact solutions and their construction}

\begin{figure}[ht]
    \centering
    \includegraphics[width=0.48\textwidth]{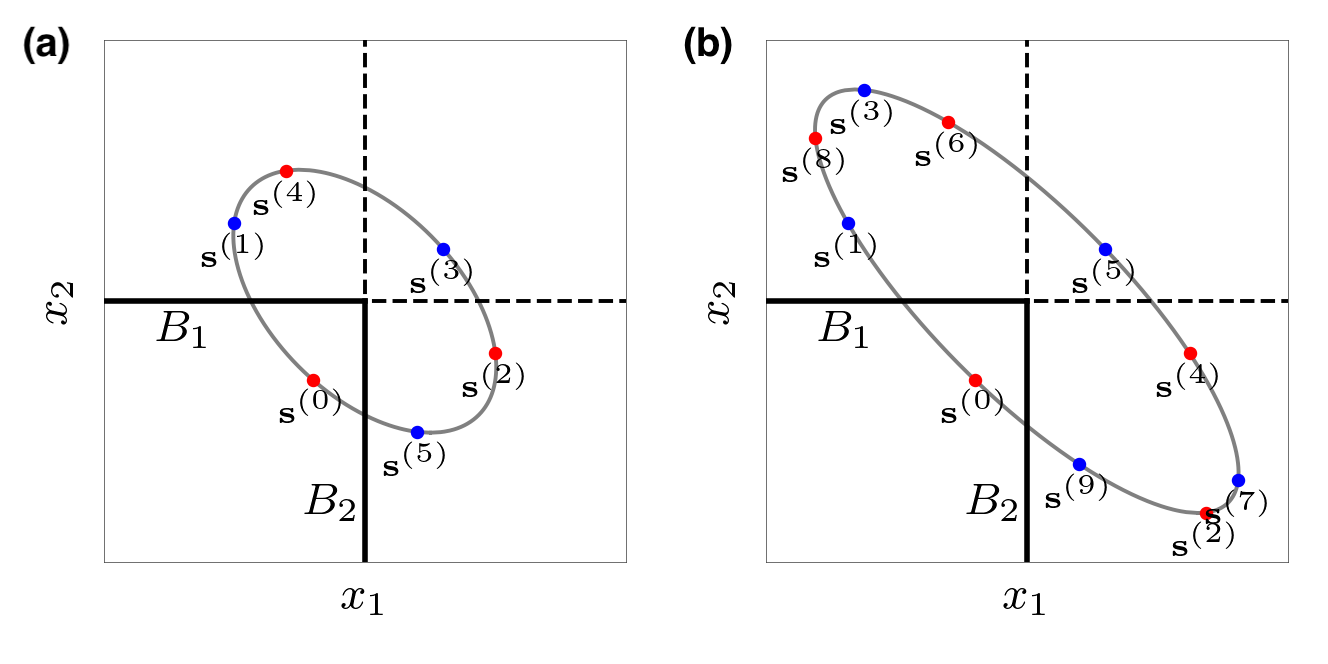}
    \caption{Illustration of full solutions in the non-whitened space. Even- and odd-numbered sources are shown as red and blue dots, respectively. (a) $k=3$ ($\rho=-\frac{1}{2}$). (b) $k=5$ ($\rho\approx -0.809$).}
    \label{fig:solution examples}
\end{figure}

Having derived the necessary and sufficient conditions on $\rho$ for $0 < |\rho| < 1$ to satisfy the initial boundary condition, we can now combine these solutions with those known for $\rho = 0$ \cite[Eq.~(A.9)]{moreno2010decision}, $\rho = 1$ \cite[Sec.5.7, Eq.~(71)]{coxmiller1965}, and $\rho = -1$ \cite[Sec.5.7, Eq.~(78)]{coxmiller1965}.
Overall, that leads to the the necessary and sufficient conditions on $\rho$ for existence of a exact solution to be given by
\begin{equation}
    \rho =  - \cos \left( \frac{\pi}{k} \right) , \quad k \in \mathbb{Z}^+ \cup \{+\infty\}.
    \label{eq:solvable rho}
\end{equation}
Where $k=\infty$ corresponds to the case of $\rho=-1$. It is worth reiterating that this condition holds \textit{regardless of}  $\hat{\bm{s}}^{(0)}$.
Solutions, if they exist, are given by the MoI construction in Table~\ref{tbl:solution_summary}.
This construction also recovers the known solutions for $\rho=-1$ (see Appendix).
Examples for $k=3$ and $k=5$ are shown in Fig.~\ref{fig:solution examples}.

\subsection{Validation with simulations}

\begin{figure}[ht]
    \centering
    \includegraphics[width=0.48\textwidth]{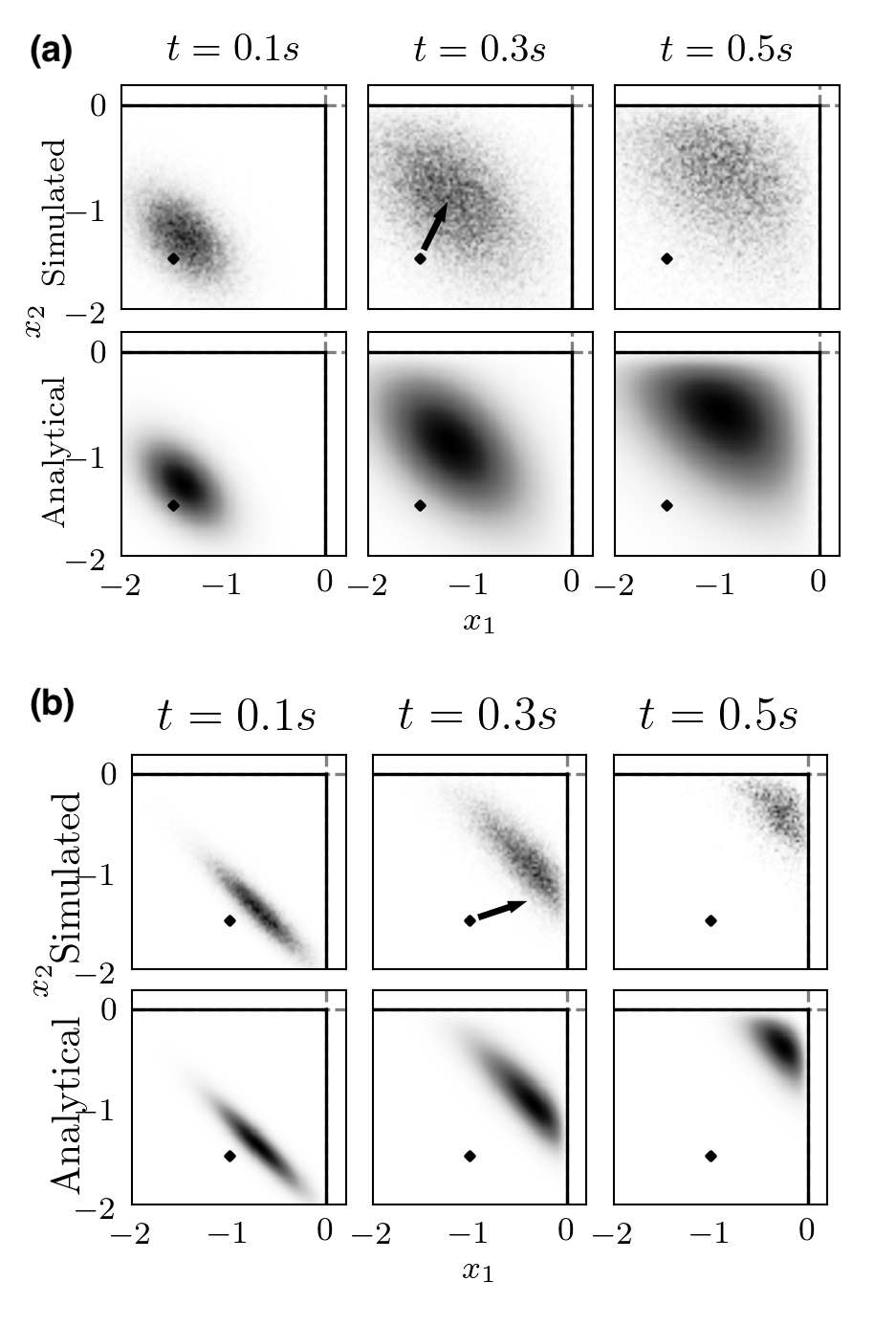}
    \caption{Illustration of PDFs at different $t$ obtained from Monte Carlo simulations (top row) and closed-form expressions (bottom row). The diamond marks $\bm{s}^{(0)}$ and the arrow indicates the direction of drift. Solid lines indicate boundaries. (a) uses $k=3$ ($\rho=-\frac{1}{2}$), $\bm{\mu} = (1,2)$, $\bm{s}^{(0)} = (-1.5,-1.5)$; (b) uses $k=8$ ($\rho\approx -0.924$), $\bm{\mu} = (3,1)$, $\bm{s}^{(0)} = (-1,-1.5)$. Note that $\bm{s}^{(0)}$ and $\bm{\mu}$ differ between (a) and (b). The analytical results were computed using the image rotation formalism. The image mapping formalism yielded, as expected, equivalent results (not shown).}
    \label{fig6}
\end{figure}
To validate the closed-form solutions we obtained, we compared them with Monte Carlo simulations based on Eq.\ (\ref{eq:langevin}), using time step-size $\delta t = 0.1ms$ and $50 000$ repetitions per figure panel. 

We first compared our closed-form expressions of $\Pdensity(\bm{x},t)$ with those obtained from simulations (Fig. \ref{fig6}). Our expressions show good agreement with results from simulation.

\begin{figure}[ht]
    \centering
    \includegraphics[width=0.48\textwidth]{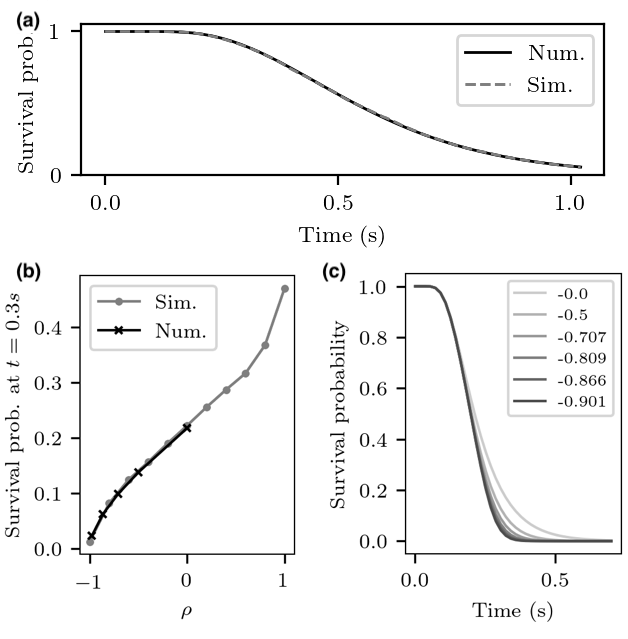}
    \caption{(a) Survival probability ($k=3$ ($\rho = -\frac{1}{2}$)). Results from numerically integrating closed-form expressions are shown in black; results from simulations are shown in gray. (b) Survival probability at a fixed time increases as a function of $\rho$. (c) Survival probability over time for some solvable $\rho$'s. In all panels, $\bm{\mu}=(2,1)$ and $\bm{s}^{(0)}=(-1.5,-1.5)$. All numerical solutions were computed using the image rotation formalism.}
    \label{fig7}
\end{figure}

A quantity of interest for diffusion processes with Dirichlet boundaries is the survival probability, defined as the probability mass within boundaries at a given time. Once these survival probabilities are known, it is easy to computer other quantities, such as the probability flux across boundaries. Using our expressions, the survival probabilities become a weighted sum of the cumulative distribution functions (CDFs) of the different images. Since each image is a scaled bivariate Gaussian distribution, evaluation of its CDFs can be carried out efficiently (we used the \texttt{stats.multivariate\_normal} object in \texttt{SciPy}). We compared the survival probability obtained from expressions to that obtained from simulations (Fig. \ref{fig7}(a)). The two again show good agreement. 

Finally, we validated by simulations that qualitative behaviors of the process vary smoothly with $\rho$, even though we could only find closed-form expressions for a limited set of $\rho$'s. To do so, we computed the survival probability from these simulations at fixed time $t = 1$ for processes with different $\rho$'s.
As shown in Fig.~\ref{fig7}(b), the survival probability varies smoothly (and for most $\rho$'s, linearly) as a function of $\rho$, and matched those found numerically for the $\rho$'s for which such numerical evaluation was possible. To further demonstrate this smoothness, we plotted the survival probability across time for various solvable $\rho$'s (Fig. \ref{fig7}(c)). The smoothness validates that we can generalizing qualitative insights from solvable $\rho$s to unsolvable ones. 

\section{Discussion}
We used the method of images (MoI) to derive a family of closed-form, analytical solutions for two-dimensional Fokker-Planck equations (FPE). The resulting solutions are unique, exact and compact. Using geometric arguments, we derived necessary and sufficient conditions for MoI solutions to exist, and validated these solutions through Monte Carlo simulations.

While we focused on two-dimensional spaces, it should be possible to generalize our approach to higher-dimensional spaces. Specifically, the analogous version of our problem in $N$-dimensions entails $N$ orthogonal, hyperplanar Dirichlet boundaries that are orthogonal to each of the axes. In such cases,  the ellipse $E(\bm{s}^{(0)})$ replaced by hyperellipsoids and lines $L_{1,2}(\bm{s}^{(0)})$ replaced by hyperplanes.
However, even for three dimensions, we would need three additional images to cancel a single image along each boundary.
For a three-dimensional, uncorrelated diffusion, for example, we require seven images, rather than the three required for an analogue two dimensional, uncorrelated diffusion.
Thus, the solution complexity will increase with the dimensionality of the space. 
A similar approach may lead to closed-form solutions for other boundary conditions, like non-orthogonal, or reflecting boundaries.

While we considered spatially homogeneous diffusion, our results can be extended to spatially heterogeneous cases as long as they can be transformed into cases we considered here via a diffeomorphism of the third quadrant. For example, we could define a new process $\bm{y}(t):y_i(t)=x_i(t)^2$ with spatially heterogeneous diffusion, since $\nabla_{\bm{y}\bm{y}}P$ depends on $\bm{y}$. Its solutions can nonetheless be found by transforming the MoI solution for $\bm{x}(t)$ with the same diffeomorphism.

The set of $\rho$'s for which we derived closed-form solutions is discrete and covers the regime of strong anti-correlations, $\rho < -1/2$ densely.
In contrast, positive correlations are not covered at all. Based on our numerical analyses, we contend that qualitative behaviors of the process are sufficiently smooth over $\rho$ that insights from solvable $\rho$'s are highly relevant in unsolvable cases as well.
Further quantitative extrapolations can be explored now by using peturbative expansions around solvable $\rho$'s, which our solutions enable.

Finally, note that showing that it is impossible to find close-form solutions for certain $\rho$'s with the MoI does not imply that there don't exist any close-form solutions for these $\rho$'s with a different from. If such solutions exist, and what form they might take, remains an open question.

\section{Acknowledgments}
This work was supported by a James S. McDonnell Foundation Scholar Award (grant \#220020462; JD), an NIMH grant (R01MH11554; JD), a MINECO (Spain) grant (BFU2017-85936-P, RMB) and a Howard Hughes Medical Institute grant (55008742, RMB).

\bibliography{citation.bib}

\onecolumngrid
\appendix
\section{Appendix}

\subsection{Properties of the whitening transformation}
\label{appendix:whitening}

\subsubsection{Source locations}
After whitening the process spatially with Eq.~(\ref{eq:whitening matrix}), the mapping between sources becomes (analogous to Eqs.~(\ref{eq:pi_plus}) and (\ref{eq:pi_minus}))
\begin{align}
    \hat{B}_1: \quad \hat{\bm{\Pi}}_1 = \bm{Q} \bm{\Pi}_1 \bm{Q}^{-1} &= \begin{pmatrix} \sqrt{1 - \rho^2} & -\rho \\ - \rho & - \sqrt{1 - \rho^2} \end{pmatrix}  , \\
    \hat{B}_2: \quad \hat{\bm{\Pi}}_2 = \bm{Q} \bm{\Pi}_2 \bm{Q}^{-1} &=  \begin{pmatrix} - \sqrt{1 - \rho^2} & -\rho \\ - \rho & \sqrt{1 - \rho^2} \end{pmatrix}.
\end{align}
Their product is
\begin{equation}
    \hat{\bm{\Pi}}_2 \hat{\bm{\Pi}}_1 = \begin{pmatrix} 2 \rho^2 - 1 & 2 \rho \sqrt{1-\rho^2} \\ - 2 \rho \sqrt{1-\rho^2} & 2 \rho^2 - 1. \end{pmatrix}=\begin{pmatrix} \cos(2\alpha) & \sin(2\alpha) \\ - \sin(2\alpha) & \cos(2\alpha) \end{pmatrix} = \bm{R}(2 \alpha),
\end{equation}
where the second equality follows from $\rho=\cos(\alpha)$ and trigonometric identities. The result is a clock-wise rotation matrix $\bm{R}(2\alpha)$ by an angle $2 \alpha$.
To find the odd-numbered sources on the circle, we will use $\hat{\bm{s}}^{(2m + 1)} = \hat{\bm{\Pi}}_1 \left(\hat{\bm{\Pi}}_2 \hat{\bm{\Pi}}_1\right)^m \hat{\bm{s}}^{(0)}$, and observe that $\hat{\bm{\Pi}}_1$ can be decomposed into
\begin{equation}
    \hat{\bm{\Pi}}_1 = \begin{pmatrix} 0 & -1 \\ -1 & 0 \end{pmatrix} \begin{pmatrix} \cos(\alpha) & \sin(\alpha) \\ -\sin(\alpha) & \cos(\alpha) \end{pmatrix} \equiv \bm{F} \left[ \bm{R}(\alpha) \right] .
\end{equation}
Here, $\bm{R}$ is again a clockwise rotation matrix, and the permutation $\bm{F}$ (see Eq.~(\ref{eq:antidiag flip})) mirrors the source across the anti-diagonal (where $x_1 = -x_2$).
Therefore, the odd-numbered source locations are given by
\begin{equation}
    \hat{\bm{s}}^{(2m + 1)} = \bm{F} \left[ \bm{R} \left( (2m + 1) \alpha \right) \hat{\bm{s}}^{(0)} \right] ,
\end{equation}
which again corresponds to an even spacing along the circle in steps of $ 2 \alpha$, but, due to the mirroring, in the opposite direction as the even-numbered sources.

\subsubsection{The avoidance region $D$}
As Eqs.~(\ref{eq:whitened_even_images}) and (\ref{eq:whitened_odd_images}) show, the sources $2m$ and $2m + 1$ are related by a one-to-one mapping.
In particular, the source $2m$ follows from a clockwise rotation by $2m \alpha$ of $\hat{\bm{s}}^{(0)}$, whereas the source $2m + 1$ follows from a similar clockwise rotation by $(2m + 1) \alpha$ of $\hat{\bm{s}}^{(0)}$, followed by a flip along the anti-diagonal.
In terms of polar angle transformation, the relations are given by 
\begin{align}
    \theta_{2m}&=\theta_0 - 2m\alpha, \\
    \theta_{2m+1}&=\frac{3}{2}\pi- \theta_{2m}+\alpha ,
    \label{eq:odd images}
\end{align}
where $\theta_{k}$ is the polar angle of $\hat{\bm{s}}^{(k)}$. 

To derive conditions for avoiding placing images in $q_{\Romanbar{III}}$, we use the above relationship to find a region $D$ such that $\theta_{2m+1}\notin q_{\Romanbar{III}} \iff \theta_{2m} \notin D$. This way, all conditions will be about even-numbered images. Expression for $D$ simply follows Eqs.~(\ref{eq:qthree_angle}) and (\ref{eq:odd images}), and result in the $D$ given by Eq.~(\ref{eq:avoidance region}) in the main text.

\subsection{Deriving the closed-form solution in the image rotation formalism}

We rely on Eqs.~(\ref{eq:whitened_even_images}) and (\ref{eq:whitened_odd_images}) to derive the closed-form solution in the image rotation formalism.
They use $2k$ images, corresponding to correlation coefficient $\rho = - \cos \left( \pi / k \right)$, resulting in the sources
\begin{equation}
    \bm{s}^{(j)} = \begin{cases}
        \bm{Q}^{-1} \bm{R} \left( j \alpha \right) \bm{Q} \bm{s}^{(0)} & j \text{ is even,} \\
        \bm{Q}^{-1} \bm{F} \bm{R} \left(j \alpha \right) \bm{Q} \bm{s}^{(0)} & j \text{ is odd,}
    \end{cases}
\end{equation}
where $\alpha = \pi (k-1) / k$,  $\bm{Q}$ and $\bm{F}$ are given by Eqs.\ (\ref{eq:whitening matrix}) and (\ref{eq:antidiag flip}), respectively, and $\bm{R}(j \alpha)$ is a 2D clockwise rotation matrix of angle $j \alpha$.
In the above, $\bm{Q}$ and $\bm{Q}^{-1}$, map into and out of the whitened space, respectively.

For even-numbered images, the image mapping can be simplified to
\begin{equation}
    \bm{Q}^{-1} \bm{R} \left( j \alpha \right) \bm{Q} = \begin{pmatrix}
        \cos \left( j \alpha \right) - \frac{\rho \sin \left( j \alpha \right)}{\sqrt{1 - \rho^2}} & \frac{\sin \left( j \alpha \right)}{\sqrt{1 - \rho^2}} \\
        - \frac{\sin \left( j \alpha \right)}{\sqrt{1 - \rho^2}} & \cos \left( j \alpha \right) + \frac{\rho \sin \left( j \alpha \right)}{\sqrt{1 - \rho^2}}
    \end{pmatrix} .
\end{equation}
Substituting $\rho = - \cos \left( \pi / k \right)$ results, after some simplification, in
\begin{equation}
    \bm{Q}^{-1} \bm{R} \left( j \alpha \right) \bm{Q} = \frac{1}{\sin \left( \frac{\pi}{k} \right)} \begin{pmatrix}
        \sin \left( j \alpha + \frac{\pi}{k} \right) & \sin \left( j \alpha \right) \\
        -\sin \left( j \alpha \right) & - \sin \left( j \alpha - \frac{\pi}{k} \right)
    \end{pmatrix} .
\end{equation}

For odd-numbered images, 
\begin{equation}
    \bm{Q}^{-1} \bm{F} \bm{R} \left( j \alpha \right) \bm{Q} = \begin{pmatrix}
        \frac{\sin \left( j \alpha \right)}{\sqrt{1 - \rho^2}} &
        - \cos \left( j \alpha \right) - \frac{\rho \sin \left( j \alpha \right)}{\sqrt{1 - \rho^2}} \\
        - \cos \left( j \alpha \right) + \frac{\rho \sin \left( j \alpha \right)}{\sqrt{1 - \rho^2}} &
        - \frac{\sin \left( j \alpha \right)}{\sqrt{1 - \rho^2}}
    \end{pmatrix} .
\end{equation}
Substituting again $\rho = - \cos \left( \pi / k \right)$ results, after some simplification, in
\begin{equation}
    \bm{Q}^{-1} \bm{F} \bm{R} \left( j \alpha \right) \bm{Q} = \frac{1}{\sin \left( \frac{\pi}{k} \right)} \begin{pmatrix}
        \sin \left( j \alpha \right) & \sin \left( j \alpha -  \frac{\pi}{k} \right) \\
        - \sin \left( j \alpha + \frac{\pi}{k} \right) & -\sin \left( j \alpha \right)
    \end{pmatrix} .
\end{equation}

\subsection{Recovering known solution for $\rho=-1$}

For $\rho = -1$, the process is a one-dimensional drift-diffusion process between two absorbing boundaries. 
The solution to this problem is provided in \cite[][Sec.~5.7, Eq.~(78)]{coxmiller1965} by using the MoI with an infinite number of images.
They denote the one-dimensional drift by $\mu$, the diffusion variance by $\sigma^2$, assume boundaries at $x = a$ and $x = -b$, and initial condition $p \left(x,0\right)=\delta\left(x\right)$.
Under these circumstances, they show the solution to be given by
\begin{equation}
p\left(x,t\right)=\frac{1}{\sigma\sqrt{2\pi t}}\sum_{k=-\infty}^{\infty}\left[\exp\left(\frac{\mu x_{k}'}{\sigma^{2}}\right)\exp\left(-\frac{\left(x-x_{k}'-\mu t\right)^{2}}{2\sigma^{2}t}\right)-\exp\left(\frac{\mu x_{k}''}{\sigma^{2}}\right)\exp\left(-\frac{\left(x-x_{k}''-\mu t\right)^{2}}{2\sigma^{2}t}\right)\right],
\end{equation}
where $x_{k}'=2k\left(a+b\right)$ and $x_{k}''=\left(2-2k\right)a-2kb$ are locations of image sources for $k=0, \pm 1, \pm 2, \dots$, and $a,b>0$ are the distances from the origin to the two boundaries.

Our formalism can also recover the infinite-image solution for $\rho=-1$. 
In this case,
\begin{equation}
    \bm{\Pi}_1 = \begin{pmatrix} 1 & 2 \\ 0 & -1 \end{pmatrix}, \qquad
    \bm{\Pi}_2 = \begin{pmatrix} -1 & 0 \\ 2 & 1 \end{pmatrix} ,
\end{equation}
placing sources on the line defined by $s_1^{(i)}+s_2^{(i)}=s_1^{(0)}+s_2^{(0)}$.
More generally, $\bm{\Pi}_1 \bm{s} = (s_1 + 2 s_2, -s_2)^T$ and $\bm{\Pi}_2 \bm{s} = (-s_1, 2s_1 + s_2)^T$.

To relate this to the above infinite-image expression, note that the drift-diffusion process is now restricted to the line $Z= \left\{  \bm{x} \big| x_{1}+x_{2} = c \right\}$, where $c = s_1^{(0)} + s_2^{(0)}$.
Let $x_Z$ denote how far we move along this line from $\bm{s}^{(0)}$ in the $(1,-1)^T$ direction, such that, for a given $x_Z$, the two-dimensional coordinates are $\bm{x} = \bm{s}^{(0)} + x_Z (1, -1)^T / \sqrt{2}$.
The mapping from $\bm{x}$ to $x_Z$ is thus given by $x_Z(\bm{x}) = \sqrt{2} \left( x_1 - s_1^{(0)} \right) = \sqrt{2} \left( s_2^{(0)} - x_2 \right)$.
This implies that the distances at which the line $Z$ intersects the boundaries, $B_1$ and $B_2$ are at
\begin{equation}
    a = x_Z \left( \left( 0, s_1^{(0)} + s_2^{(0)} \right)^T \right) = - \sqrt{2} s_1^{(0)}, \qquad \text{ and } \qquad
    b = - x_Z \left( \left( s_1^{(0)} + s_2^{(0)}, 0 \right)^T \right) = - \sqrt{2} s_2^{(0)} ,
\end{equation}
respectively.
Furthermore, it is easy to verify that
\begin{align}
    x_Z \left( \bm{\Pi}_1 \bm{x} \right) &= - 2 b - x_Z \left( \bm{x} \right) , \\
    x_Z \left( \bm{\Pi}_2 \bm{x} \right) &= 2 a - x_Z \left( \bm{x} \right) .
\end{align}
Then $x_Z \left( \bm{s}^{(0)} \right) = x_0'$, and
\begin{align}
    x_Z \left( \bm{\Pi}_1 \bm{s}^{(0)} \right), x_Z \left( \bm{\Pi}_2 \bm{\Pi}_1 \bm{s}^{(0)} \right), x_Z \left( \bm{\Pi}_1 \bm{\Pi}_2 \bm{\Pi}_1 \bm{s}^{(0)} \right), x_Z \left( \bm{\Pi}_2 \bm{\Pi}_1 \bm{\Pi}_2 \bm{\Pi}_1 \bm{s}^{(0)} \right), \dots &= -2b, 2a + 2b, -2a - 4b, 4a + 4b, \dots \nonumber \\
    &= x_1'', x_1', x_2'', x_2', \dots, \\
    x_Z \left( \bm{\Pi}_2 \bm{s}^{(0)} \right), x_Z \left( \bm{\Pi}_1 \bm{\Pi}_2 \bm{s}^{(0)} \right), x_Z \left( \bm{\Pi}_2 \bm{\Pi}_1 \bm{\Pi}_2 \bm{s}^{(0)} \right), x_Z \left( \bm{\Pi}_1 \bm{\Pi}_2 \bm{\Pi}_1 \bm{\Pi}_2 \bm{s}^{(0)} \right), \dots &= 2a, -2a-2b, 4a+2b, -4a-4b, \dots, \nonumber \\
    &=x_0'', x_{-1}', x_{-1}'', x_{-2}', \dots
\end{align}
which corresponds to the image sequence of the above solution.

\end{document}